\documentclass[a4paper,twoside,11pt]{article}
\usepackage{a4wide, fancyhdr, amsmath, amssymb, mathtools, yfonts}
\usepackage{graphicx}
\usepackage{tikz}
\usepackage[all]{xy}
\usepackage[utf8]{inputenc}
\usepackage{amsthm}
\usepackage[english]{babel}
\usepackage{chngcntr}
\usepackage{ifthen}
\usepackage{calc}
\usepackage{hyperref}
\usepackage[OT2,T1]{fontenc}
\usepackage{authblk}
\numberwithin{equation}{section}
\DeclareSymbolFont{cyrletters}{OT2}{wncyr}{m}{n}
\DeclareMathSymbol{\Sha}{\mathalpha}{cyrletters}{"58}

\providecommand{\keywords}[1]{{2010}\textit{ Mathematics Subject Classification. }#1}

\newcommand{\ZZ}{\mathbb{Z}}
\newcommand{\QQ}{\mathbb{Q}}

\newcommand{\OO}{\mathcal{O}}
\newcommand{\CL}{\mathrm{Cl}}

\newcommand\ZZT{\mathbb{Z}[\sqrt{2}]}
\newcommand\ZZK{\mathbb{Z}[\zeta_8]}

\newcommand\QQT{\mathbb{Q}(\sqrt{2})}

\newcommand\mm{\mathfrak{m}}

\newcommand\nn{\mathfrak{n}}
\newcommand\dd{\mathfrak{d}}
\newcommand\pp{\mathfrak{p}}

\newcommand\pt{\mathfrak{t}}

\newcommand\bbb{\mathfrak{b}}

\newcommand\Gal{\mathrm{Gal}}
\newcommand\Norm{\mathfrak{N}}

\newcommand\ve{\varepsilon}

\newtheorem{theorem}{Theorem}
\newtheorem{lemma}{Lemma}[section]
\newtheorem{prop}[lemma]{Proposition}
\newtheorem{corollary}[theorem]{Corollary}

\title{On the equations $x^2-2py^2 = -1, \pm 2$}
\author{Djordjo Z. Milovic\thanks{Gower Street, London, WC1E 6BT, United Kingdom, djordjo.milovic@ucl.ac.uk}}
\affil{Department of Mathematics, University College London}

\date{\today}

\begin{document}

\maketitle

\begin{abstract} 
Let $E\in\{-1, \pm 2\}$. We improve on the upper and lower densities of primes $p$ such that the equation $x^2-2py^2=E$ is solvable for $x, y\in \ZZ$. We prove that the natural density of primes $p$ such that the narrow class group of the real quadratic number field $\QQ(\sqrt{2p})$ has an element of order $16$ is equal to $\frac{1}{64}$. We give an application of our results to the distribution of Hasse's unit index for the CM-fields $\QQ(\sqrt{2p}, \sqrt{-1})$. Our results are consequences of a twisted joint distribution result for the $16$-ranks of class groups of $\QQ(\sqrt{-p})$ and $\QQ(\sqrt{-2p})$ as $p$ varies.
\end{abstract}

\section{Introduction}\label{Mainresults}
Let $p$ denote an odd prime number. Gauss's genus theory implies that the $2$-torsion subgroup of the narrow class group $\CL^+(2p)$ of the real quadratic number field $\QQ(\sqrt{2p})$ is isomorphic to the group of two elements, and that it is generated by the classes of the ideals $\pt = (2, \sqrt{2p})$ and $\pp = (p, \sqrt{2p})$ in $\ZZ[\sqrt{2p}]$ \cite[Lemma 9.8(a), p.78]{Ste1}. Hence exactly one of the ideals $\pt$, $\pp$, and $\pt\pp = (\sqrt{2p})$ is principal in the narrow sense, and the remaining two are in the class of order $2$ in $\CL^+(2p)$. Let $E_p$ denote the unique integer in the set $\{-1, \pm 2\}$ such that the equation
\begin{equation}\label{gPE}
x^2 - 2py^2 = E_p
\end{equation} 
has a solution with $x, y\in \ZZ$; in other words, if we denote the class of an ideal $\nn\subset \ZZ[\sqrt{2p}]$ in $\CL^+(2p)$ by $[\nn]$, then
$$
E_p = 
\begin{cases}
-1 & \text{if }[\pt\pp] = 1, \\
2 & \text{if }[\pt] = 1, \\
-2 & \text{if }[\pp] = 1.
\end{cases}
$$
Stevenhagen conjectures that each of the three cases above occurs equally often \cite[p.127]{SteConj}. More precisely, let $E\in\{-1, \pm 2\}$, let $X\geq 3$ be a real number, and let $\delta(X)$ denote the proportion 
\begin{equation}\label{defdE}
\delta(X; E) = \frac{|\{p\leq X: E_p = E\}|}{|\{p\leq X\}|}.
\end{equation}
Then Stevenhagen conjectures that the limit $\lim_{X\rightarrow\infty}\delta(X; E)$ exists and is equal to $\frac{1}{3}$. Classical results of R\'{e}dei \cite{Redei2} and Scholz \cite{Scholz} from the 1930's can be used to deduce that 
\begin{equation}\label{classical}
\frac{5}{16}\leq \liminf_{X\rightarrow\infty}\delta(X; E)\leq \limsup_{X\rightarrow\infty}\delta(X; E) \leq \frac{3}{8},
\end{equation}
and Koymans and the author \cite{KM2} recently proved that
$$
\limsup_{X\rightarrow\infty}\delta(X; -1) \leq \frac{11}{32},
$$
conditional on Conjecture~$C_n$ with $n = 8$ from \cite{FIMR}; thanks to the recent work of Koymans \cite{Koymans}, one can remove this conditionality on Conjecture~$C_n$. In this paper, we prove the same upper bounds for $E = \pm 2$ as well as improve the lower bounds in all cases. Our results are unconditional.
\begin{theorem}\label{THM1}
Let $E\in \{-1, \pm 2\}$, let $X\geq 3$ be a real number, and define $\delta(X; E)$ as in \eqref{defdE}. Then
$$
\frac{21}{64}\leq \liminf_{X\rightarrow\infty}\delta(X; E)\leq \limsup_{X\rightarrow\infty}\delta(X; E) \leq \frac{11}{32}.
$$ 
\end{theorem}
Furthermore, we will see that the ``remaining'' primes, which form a set of natural density $1-3\cdot \frac{21}{64} = \frac{1}{64}$, are exactly those for which $\CL^+(2p)$ has an element of order $16$. This yields the first density result on the $16$-rank in any family of \textit{real} quadratic fields parametrized by one prime number. Recall that the $2$-part of $\CL^+(2p)$ is cyclic and set $h^+(-2p) = |\CL^+(2p)|$.
\begin{theorem}\label{THM2}
We have
$$
\lim_{X\rightarrow\infty}\frac{|\{p\leq X: h^+(2p)\equiv 0\bmod 16\}|}{|\{p\leq X\}|}= \frac{1}{64}.
$$ 
\end{theorem}
Finally, we also have an application to the distribution of Hasse's unit index for the biquadratic fields $\QQ(\sqrt{2p}, \sqrt{-1})$. Let $U_F$ denote the unit group of the ring of integers of an algebraic number field $F$. Let $K$ be a CM-field with maximal real subfield $K^+$ and let $W_K$ denote the subgroup of $U_K$ consisting of units of finite order. Then \textit{Hasse's unit index} $Q(K)$ of $K$ is defined to be the index $[U_K:U_{K^+}W_K]$; it is always equal to $1$ or $2$. If $K = \QQ(\sqrt{p}, \sqrt{-1})$, then $Q(K) = 1$, so the simplest case from the standpoint of arithmetic statistics is when $K = \QQ(\sqrt{2p}, \sqrt{-1})$. In that case, $Q(K) = 1$ if and only if a fundamental unit $\epsilon_{2p}$ of $\ZZ[\sqrt{2p}]$ and $\sqrt{-1}$ generate the full unit group $U_K$. This occurs if and only if the ideal $\pt$ from above is \textit{not} principal in the ordinary sense \cite[Theorem 1.ii.2.]{Lemm}, i.e., if and only if $E_p \neq \pm 2$. Hence we can deduce the following corollary of Theorem~\ref{THM1}.
\begin{corollary}
Let $K = \QQ(\sqrt{2p}, \sqrt{-1})$, let $U_K$ denote the group of its integral units, and let $\epsilon_{2p}$ denote a fundamental unit of $\ZZ[\sqrt{2p}]$. For a real number $X\geq 3$, let 
$$
\delta_{H}(X) =  \frac{|\{p\leq X: U_K\text{ is generated by }\epsilon_{2p}\text{ and }\sqrt{-1}\}|}{|\{p\leq X\}|}.
$$
Then
$$
\frac{21}{64}\leq \liminf_{X\rightarrow\infty}\delta_H(X)\leq \limsup_{X\rightarrow\infty}\delta_H(X) \leq \frac{11}{32}.
$$
\end{corollary}  

\section{Main strategy}
Our results are possible in large part thanks to the work of Kaplan and Williams \cite{KW84} and Leonard and Williams \cite{LW82}. Kaplan and Williams relate the existence of an element of order $16$ in the narrow class group of the \textit{real} quadratic field $\QQ(\sqrt{2p})$ to the existence of elements of order $16$ in the class groups of the \textit{imaginary} quadratic fields $\QQ(\sqrt{-p})$ and $\QQ(\sqrt{-2p})$. A similar type of ``reflection principle'' for the $16$-rank was also later proved by Stevenhagen \cite{Ste2}, but the results from \cite{KW84} appear to be more suitable for our purposes and our analytic techniques. 

Let $h(-p)$ and $h(-2p)$ denote the class numbers of $\QQ(\sqrt{-p})$ and $\QQ(\sqrt{-2p})$, respectively. Note that Gauss's genus theory implies that the $2$-parts of the narrow class groups of $\QQ(\sqrt{-p})$, $\QQ(\sqrt{-2p})$, and $\QQ(\sqrt{2p})$ are cyclic and hence determined by the highest power of $2$ dividing $h(-p)$, $h(-2p)$, and $h^+(2p)$, respectively. The following facts can be found in or readily deduced from \cite{KW84, Ste2}. We have
\begin{align*}
h^+(2p)\equiv 0\bmod 8 & \Leftrightarrow h(-p) \equiv h(-2p)\equiv 0\bmod 8 \\
& \Leftrightarrow p\text{ splits completely in }\QQ(\zeta_{16}, \sqrt[4]{2}).
\end{align*}
Here $\zeta_{n}$ denotes a primitive $n$th root of unity. Let $F$ be any one of the three quadratic fields $\QQ(\sqrt{-1})$, $\QQ(\sqrt{-2})$, and $\QQ(\sqrt{2})$. It is not hard to check that in each case 
$$
\Gal(\QQ(\zeta_{16}, \sqrt[4]{2}))/F)\cong \ZZ/2\times\ZZ/2\times \ZZ/2.
$$
As this group is abelian, class field theory implies that if a prime $p$ that splits in $F/\QQ$, then the splitting type of $p$ in $\QQ(\zeta_{16}, \sqrt[4]{2})/\QQ$ can be detected via congruence conditions on a prime $\pi$ in $F$ lying above $p$. Concretely, $p$ splits completely in $\QQ(\zeta_{16}, \sqrt[4]{2})/\QQ$ if and only if there exist integers $a$, $b$, $c$, $d$, $u$, and $v$ such that 
\begin{equation}\label{z1}
p = a^2+b^2 = c^2+2d^2 = u^2-2v^2
\end{equation}
and
\begin{equation}\label{z2}
a\equiv 1\bmod 8,\ b\equiv 0\bmod 8,\ c\equiv 1\bmod 8,\ d\equiv 0\bmod 4,\ u\equiv 1\bmod 8,\ v\equiv 0\bmod 4.
\end{equation}
For such a prime $p$, define
\begin{equation}\label{defapbp}
\alpha_p = (-1)^{h(-2p)/8}\quad\text{and}\quad\beta_p = (-1)^{(a-1+b+2d+h(-p)+h(-2p))/8}.
\end{equation}
By studying the number of quadratic residues modulo $p$ that are less than $p/8$, Kaplan and Williams \cite[Theorem, p.26]{KW84} proved that 
\begin{align}\label{21}
&\alpha_p = \beta_p = 1 & \Longrightarrow & h^+(2p)\equiv 0\bmod 16,\\ \label{22}
&\alpha_p = 1, \beta_p = -1 & \Longrightarrow & h^+(2p)\equiv 8\bmod 16,\ E_p = -2, \\ \label{23}
&\alpha_p = -1, \beta_p = 1 & \Longrightarrow & h^+(2p)\equiv 8\bmod 16,\ E_p = +2, \\ \label{24}
&\alpha_p = \beta_p = -1 & \Longrightarrow & h^+(2p)\equiv 8\bmod 16,\ E_p = -1.
\end{align}
We will prove that each of the four possibilities \eqref{21}-\eqref{24} occurs equally often, and this will imply both Theorem~\ref{THM1} and Theorem~\ref{THM2}. Indeed, the natural density of primes $p$ that split completely in $\QQ(\zeta_{16}, \sqrt[4]{2})/\QQ$ is equal to $\frac{1}{16}$, by the Chebotarev Density Theorem; hence we will prove that the natural density of primes satisfying each of the four possibilities above is equal to $\frac{1}{64}$. This immediately implies Theorem~\ref{THM2}. We note that the classical bounds \eqref{classical} derive from primes $p$ that do \textit{not} split completely in $\QQ(\zeta_{16}, \sqrt[4]{2})/\QQ$ (this is a set of primes of natural density $\frac{15}{16}$), and the improvements in Theorem~\ref{THM1}, say for $E = -2$, come from adding to the lower bound in \eqref{classical} the fraction of primes satisfying \eqref{22} and subtracting from the upper bound in \eqref{classical} the fraction of primes satisfying \eqref{23} or \eqref{24}.

To prove that each of the four possibilities \eqref{21}-\eqref{24} occurs equally often, we restrict to the set of primes $p$ that split completely in $\QQ(\zeta_{16}, \sqrt[4]{2})/\QQ$ and consider indicator functions of the type
\begin{equation}\label{indicator}
\frac{1}{4}\left(1+\alpha_p+\beta_p+\alpha_p\beta_p\right)=
\begin{cases}
1 & \text{if }\alpha_p=\beta_p=1, \\
0 & \text{otherwise},
\end{cases}
\end{equation}
It then suffices to prove that each of the three sums 
$$
\sum_{p\leq X}^{\ast}\alpha_p, \quad \sum_{p\leq X}^{\ast}\beta_p,\quad\sum_{p\leq X}^{\ast}\alpha_p\beta_p
$$
is $o(X/\log X)$ as $X\rightarrow \infty$; here $\ast$ denotes the restriction to primes that split completely in $\QQ(\zeta_{16}, \sqrt[4]{2})/\QQ$.

The sum $\sum^{\ast} \alpha_p$ encodes the behavior of the $16$-rank in the family $\{\QQ(\sqrt{-2p}): p\equiv 1\bmod 4\}$; it was the subject of a paper of Koymans and the author \cite{KM1}, where we proved that
$$
\sum_{p\leq X}^{'}\alpha_p \ll X^{\frac{1}{3200}}.
$$
Here $'$ restricts the sum to primes $p\equiv 1\bmod 4$ such that $h(-2p)\equiv 0\bmod 8$; these are exactly the primes that split completely in $\QQ(\zeta_8, \sqrt[4]{2})$. We will show that the same result holds when we further restrict the sum to $p$ such that $h^+(2p)\equiv 0\bmod 8$, i.e., $p$ that split completely in $\QQ(\zeta_{16}, \sqrt[4]{2})$.
\begin{prop}\label{PropAp}
We have
$$
\sum_{p\leq X}^{\ast}\alpha_p \ll X^{1-\frac{1}{3200}},
$$
where $\ast$ restricts the sum to primes $p$ that split completely in $\QQ(\zeta_{16}, \sqrt[4]{2})$, and where $\alpha_p$ is defined in \eqref{defapbp}.
\end{prop}
%

The sum $\sum^{\ast} \alpha_p\beta_p$ concerns a twisted version of the $16$-rank in the family $\{\QQ(\sqrt{-p})\}$; the twist here is the factor $(-1)^{(a-1+b+2d)/8}$. This sum too was the subject of another paper of Koymans and the author \cite{KM2}, where we proved that
$$
\sum_{p\leq X}^{*}\alpha_p\beta_p \ll X^{1-\frac{\delta}{400}}.
$$
Here $\delta>0$ is a conjectural constant that appears in Conjecture $C_n$ for $n=8$ in \cite{FIMR}; see \cite[Theorem 3, p.102]{KM2} and its proof in \cite[Section 7]{KM2}. Koymans~\cite{Koymans} recently gave an unconditional proof of \cite[Theorem 2, p.102]{KM1}, i.e., a similar power-saving estimate for the sum $\sum^{\ast} (-1)^{h(-p)/8}$. We will show that the proof in \cite{Koymans} can be modified slightly to also give an unconditional estimate for $\sum^{\ast}\alpha_p\beta_p$.
\begin{prop}\label{PropApBp}
We have
$$
\sum_{p\leq X}^{\ast}\alpha_p\beta_p \ll X^{1-\frac{1}{25000}},
$$
where $\ast$ restricts the sum to primes $p$ that split completely in $\QQ(\zeta_{16}, \sqrt[4]{2})$, and where $\alpha_p$ and $\beta_p$ are defined in \eqref{defapbp}
\end{prop}
It remains to show a similar estimate for the sum $\sum^{\ast}\beta_p$, which concerns a twisted version of the \textit{joint distribution} of the $16$-ranks in the families $\{\QQ(\sqrt{-p})\}$ and $\{\QQ(\sqrt{-2p})\}$. This forms the main subject of the present paper. We will prove
\begin{prop}\label{PropBp}
We have
$$
\sum_{p\leq X}^{\ast}\beta_p \ll X^{1-\frac{1}{200}},
$$
where $\ast$ restricts the sum to primes $p$ that split completely in $\QQ(\zeta_{16}, \sqrt[4]{2})$, and where $\beta_p$ is defined in \eqref{defapbp}.
\end{prop}
Propositions~\ref{PropAp}, \ref{PropApBp}, and \ref{PropBp}, in conjunction with \eqref{indicator} and similar identities, imply Theorems~\ref{THM1} and \ref{THM2}. We start by laying the groundwork for the proof of Proposition~\ref{PropBp}; with the appropriate set-up, Propositions~\ref{PropAp} and \ref{PropApBp} will follow readily from \cite{KM1} and \cite{Koymans}, respectively. Ultimately we will see that Proposition~\ref{PropBp} is connected to the equidistribution results from \cite{Milovic2}, although this connection is far less obvious -- in fact, \cite{Milovic2} features a result for a family parametrized by primes $p\equiv -1\bmod 4$, while the sum in Proposition~\ref{PropBp} is supported on primes $p\equiv 1\bmod 4$.

\section{Algebraic Criteria}\label{SectionAC}
Leonard and Williams \cite{LW82} use the theory of binary quadratic forms, along with composition laws of Gauss and Dirichlet, to derive formulas for $(-1)^{h(-p)/8}$ and $(-1)^{h(-2p)/8}$ for $p\equiv 1\bmod 8$ in terms of prime ideals lying above $p$ in $\ZZT$. The proofs of their results rely in part on very clever manipulations of Legendre and Jacobi symbols, and they produce formulas which vaguely resemble two different types of \textit{spin symbols}, one appearing in the work of Friedlander and Iwaniec \cite[(20.1), p.1021]{FI1} and the other being the main object of study by Friedlander, Iwaniec, Mazur, and Rubin in \cite{FIMR}. In the case of a similar criterion for $(-1)^{h(-2p)/8}$ for $p\equiv -1\bmod 8$ \cite[Theorem 3, p.205]{LW82}, we translated their proof to the language of rings and ideals, and this translation revealed additional structure that allowed us to embed a Jacobi symbol appearing in this criterion into a sequence conducive to certain sieving methods; see \cite{Milovic2}. As a result, it may be interesting to translate also the proofs of the following criteria to the language of rings and ideals, but we avoid doing so since the results of \cite{Milovic2} already suffice for our applications. 

\subsection{Preliminaries}
Given an integer $n\geq 1$, let $\zeta_n$ denote a primitive $n$-th root of unity. Let $F$ be a finite Galois extension of $\QQ$ containing $\zeta_n$, and let $\OO_F$ denote the ring of integers of $F$. Let $\Norm_{F/\QQ}$ denote the norm map from $F$ to $\QQ$. Given $\alpha\in\OO_F$ and a prime ideal $\pp$ in $\OO_F$ coprime to $n$, the $n$-th power residue symbol $(\alpha/\pp)_{F, n}$ is defined to be the unique element of $\{0, 1, \zeta_n, \zeta_n^2, \ldots, \zeta_n^{n-1}\}$ such that
$$
\left(\frac{\alpha}{\pp}\right)_{F, n} \equiv \alpha^{\frac{\Norm_{F/\QQ}(\pp)-1}{n}}\bmod \pp.
$$
It is evident from this definition that $(\alpha/\pp)_{F, n}$ depends only on the congruence class of $\alpha$ modulo $\pp$; that $(\alpha_1\alpha_2/\pp)_{F, n} = (\alpha_1/\pp)_{F, n}(\alpha_2/\pp)_{F, n}$ for all $\alpha_1, \alpha_2\in\OO_F$; and that $\alpha$ is an $n$-th power modulo $\pp$ if and only if $(\alpha/\pp)_{F, n} = 1$. Moreover, if $n$ is even, $\alpha$ is an $(n/2)$-th power but not an $n$-th power modulo $\pp$ if and only if $(\alpha/\pp)_{F, n} = -1$. For an ideal $\bbb$ in $\OO_F$ coprime to $n$, we set $(\alpha/\bbb)_{F, n} = \prod_{\pp^{e_{\pp}}\| \bbb}(\alpha/ \pp)_{F, n}^{e_{\pp}}$, where $\pp^{e_{\pp}}$ is the exact power of $\pp$ dividing $\bbb$. For an odd element $\beta\in\OO_F$, we set $(\alpha/\beta)_{F, n} = (\alpha/\beta\OO_F)_{F, n}$. If $F = \QQ$, $n = 2$, and $\beta\in\ZZ$ is positive, the symbol $(\cdot/ \beta)_{\QQ, 2}$ coincides with the usual Jacobi symbol, and so we suppress the subscripts $\QQ, 2$ and simply write $(\cdot/ \beta)$.

\subsection{Symbols over $\ZZ$}
Let $\ve = 1+\sqrt{2}$, and note that $\ve$ is a unit of infinite order in $\ZZT$. For the remainder of Section~\ref{SectionAC}, let $p\equiv 1\bmod 8$ be a prime number. Since $p$ splits in the unique factorization domain $\ZZ[\sqrt{2}]$, there exist rational integers $u$ and $v$ such that
\begin{equation}\label{defuv}
p = (u+v\sqrt{2})(u-v\sqrt{2}) = u^2 - 2v^2,\quad u, v>0.
\end{equation}
Note that $u$ and $v$ must be odd and even, respectively. Moreover, after multiplying $u+v\sqrt{2}$ by $\ve^2$ if necessary, we can choose $u$ in \eqref{defuv} so that 
\begin{equation}\label{umod4}
u\equiv 1 \bmod 4.
\end{equation}
Let $g$ and $h$ be positive rational integers such that
$$
h+g\sqrt{2} = (u+v\sqrt{2})\cdot\ve,
$$
so that $p = 2g^2 - h^2$. Now assume also that $h(-p)\equiv h(-2p)\equiv 0\bmod 8$ so that $p$ splits completely in $\QQ(\zeta_{16}, \sqrt[4]{2})$ and $h^+(2p)\equiv 0\bmod 8$. Then Leonard and Williams show that
\begin{equation}\label{crit8rank1}
\left(\frac{u}{p} \right) = \left(\frac{-2}{u} \right) = \left(\frac{g}{p} \right) = \left(\frac{-1}{g} \right) = 1;
\end{equation}
in particular, $u\equiv 1\bmod 8$ and $g\equiv 1\bmod 4$. As $p\equiv 1\bmod 16$, we see that $v\equiv 0\bmod 4$. Recall also \eqref{z1} and \eqref{z2} above.

For an integer $n$ satisfying $(n/p) = 1$, set
$$
\left[\frac{n}{p} \right]_4 = 
\begin{cases}
1 & \text{if }n\text{ is a fourth power residue modulo }p \\
-1 & \text{otherwise.}
\end{cases}
$$
Two of the main results in \cite{LW82} are then as follows. Suppose $u$ satisfies \eqref{defuv}-\eqref{umod4}. First, \cite[Theorem 2, p.204]{LW82} implies that
$$
(-1)^{h(-2p)/8} = \left[\frac{u}{p} \right]_4.
$$
Second, \cite[Theorem 1, p.201]{LW82} implies that
$$
(-1)^{h(-p)/8} = \left[\frac{g}{p} \right]_4\left(\frac{2h}{g} \right).
$$
Hence,
\begin{equation}\label{overZ}
(-1)^{(h(-p)+h(-2p))/8} = \left[\frac{u}{p} \right]_4\left[\frac{g}{p} \right]_4\left(\frac{2h}{g} \right).
\end{equation}
\subsection{From $\ZZ$ to $\ZZK$}
To fully exploit the multiplicative properties underlying the symbol $[\cdot / p]_4$, we will rewrite the above criterion in a field containing $\zeta_4 = \sqrt{-1}$, a primitive fourth root of unity. Moreover, as the criterion naturally depends on the splitting of $p$ in $\QQT$, we will work over the field 
$$
K= \QQ(\sqrt{-1}, \sqrt{2}) = \QQ(\zeta_8).
$$
Let $\ve = 1+\sqrt{2}$. Note that the ring of integers of $K$ is $\ZZK$, that $\ZZK$ is a principal ideal domain, and that its group of units is generated by $\zeta_8$ and $\ve$.

Since the prime $p\equiv 1\bmod 8$ splits completely in $K$, we can choose a prime element $\varpi\in\ZZK$ such that $\Norm_{K/\QQ}(\varpi) = p$. Let $\sigma$ and $\tau$ be the non-trivial elements of $\Gal(K/\QQ)\cong V_4$ fixing $\QQ(\sqrt{2})$ and $\QQ(\sqrt{-1})$, respectively. Then
$$
u + v\sqrt{2} = \varpi\sigma(\varpi), \quad h + g\sqrt{2} = \ve\varpi\sigma(\varpi),
$$ 
and so, as $\tau$ acts non-trivially on $\QQT$, we have
$$
u = \frac{1}{2}\left(\varpi\sigma(\varpi) + \tau(\varpi)\tau\sigma(\varpi)\right)
$$
and
$$
g = \frac{1}{2\sqrt{2}}\left(\ve\varpi\sigma(\varpi) + \ve^{-1}\tau(\varpi)\tau\sigma(\varpi)\right).
$$
Now let $\varpi$ be an element of norm $p$ as above. As $\varpi\ZZK$ is a prime ideal of degree one, the inclusion of rings $\ZZ\hookrightarrow \ZZK$ induces an isomorphism of fields $\ZZ/p\ZZ\cong \ZZK/\varpi\ZZK$, so an element $n\in\ZZ$ is a square (resp.\ a fourth power) modulo $p$ if and only if $n$, viewed as an element of $\ZZK$, is a square (resp.\ a fourth power) modulo $\varpi\ZZK$.  

Hence, as $4 = \sqrt{2}^4$ is a fourth power in $\ZZK$, we have
$$
\left[\frac{u}{p}\right]_4 = \left(\frac{u}{\varpi}\right)_{K, 4} = \left(\frac{2\tau(\varpi)\sigma\tau(\varpi)}{\varpi}\right)_{K, 4}
$$  
and
$$
\left[\frac{g}{p}\right]_4 = \left(\frac{g}{\varpi}\right)_{K, 4} = \left(\frac{\sqrt{2}\ve^{-1}\tau(\varpi)\sigma\tau(\varpi)}{\varpi}\right)_{K, 4}.
$$  
Combining the above formulas, and noting that
$$
\left(\frac{\tau(\varpi)\sigma\tau(\varpi)}{\varpi}\right)_{K, 4}^2 = \left(\frac{\tau(\varpi)\sigma\tau(\varpi)}{\varpi}\right)_{K, 2},
$$
one can rewrite \eqref{overZ} as 
\begin{equation}\label{overZ8}
(-1)^{(h(-p)+h(-2p))/8}=\left(\frac{2\sqrt{2}\ve^{-1}}{\varpi}\right)_{K, 4}\left(\frac{\tau(\varpi)\sigma\tau(\varpi)}{\varpi}\right)_{K, 2}\left(\frac{2h}{g}\right).
\end{equation}

\subsection{The first factor}
We have
$$
\left(\frac{2\sqrt{2}\ve^{-1}}{\varpi}\right)_{K, 4} = \left(\frac{2(2-\sqrt{2})}{\varpi}\right)_{K, 4} = \left(\frac{2}{\varpi}\right)_{K, 4}\left(\frac{e_2}{\varpi}\right)_{K, 4},
$$
where $e_2 = 2-\sqrt{2}$ as in \cite[(5), p.24]{KW84}. As $\varpi$ lies above a prime $p$ that splits completely in $\QQ(\zeta_{16}, \sqrt[4]{2})$, we see that $2$ is a fourth power modulo $p$. Thus
$$
\left(\frac{2}{\varpi}\right)_{K, 4} = \left[\frac{2}{p}\right]_4 = 1.
$$ 
Moreover, Kaplan and Williams already computed in \cite[p.25]{KW84} that
$$
\left(\frac{e_2}{\varpi}\right)_{K, 4} = \left[\frac{e_2}{p} \right]_4 = (-1)^{(b+2d)/8}.
$$
Hence 
\begin{equation}\label{firstfactor}
\left(\frac{2\sqrt{2}\ve^{-1}}{\varpi}\right)_{K, 4} = (-1)^{(b+2d)/8}.
\end{equation}

\subsection{The middle factor}
We now deal with the middle factor on the right-hand-side of \eqref{overZ8}. As $\varpi$ is a prime of degree one in $\ZZK$, setting $\pi = \varpi\sigma(\varpi) = \Norm_{\QQ(\zeta_8)/\QQT}(\varpi)$, the inclusion of rings $\ZZT\hookrightarrow \ZZK$ induces an isomorphism of fields $\ZZT/\pi\ZZT\cong \ZZK/\varpi\ZZK$. Letting $\overline{\tau}$ denote the restriction of $\tau$ to $\QQT$, we note that $\tau(\varpi)\sigma\tau(\varpi) = \tau(\varpi\sigma(\varpi)) = \overline{\tau}(\pi)\in\ZZT$, and so
$$
\left(\frac{\tau(\varpi)\sigma\tau(\varpi)}{\varpi}\right)_{K, 2} = \left(\frac{\overline{\tau}(\pi)}{\varpi}\right)_{K, 2} =  \left(\frac{\overline{\tau}(\pi)}{\pi}\right)_{\QQT, 2}
$$
Writing $\pi = u+v\sqrt{2}$ as above, we see that
$$
\left(\frac{\overline{\tau}(\pi)}{\pi}\right)_{\QQT, 2} = \left(\frac{\overline{\tau}(\pi)+\pi}{\pi}\right)_{\QQT, 2} = \left(\frac{2u}{\pi}\right)_{\QQT, 2}.
$$ 
Again, as $\pi$ is a prime of degree one and $2u\in\ZZ$, we can use the canonical isomorphism $\ZZ/p\ZZ\cong \ZZT/\pi\ZZT$ to write
$$
\left(\frac{2u}{\pi}\right)_{\QQT, 2} = \left(\frac{2u}{p}\right).
$$
By \eqref{crit8rank1}, we see that the Legendre symbol above is equal to $1$, and so the middle factor in \eqref{overZ8} is trivial, i.e.,
\begin{equation}\label{middlefactor}
\left(\frac{\tau(\varpi)\sigma\tau(\varpi)}{\varpi}\right)_{K, 2} = 1.
\end{equation}

\subsection{The last factor: essential spin}
Next, we deal with the last factor on the right-hand-side of \eqref{overZ8}. This factor is essential in ensuring that $\beta_p$ can be written as a genuine spin symbol. We now relate it to the spin symbol
$$
[u+v\sqrt{2}] = \left(\frac{v}{u}\right)
$$
appearing in \cite{Milovic2}. Recall that $u\equiv 1\bmod 8$, that $v\equiv 0\bmod 4$, and that $u$ and $v$ are positive. Writing $e_v$ for the highest power of $2$ dividing $v$, so that $v' = v2^{-e_v}$ is odd, we have
$$
\left(\frac{2h}{g}\right) = \left(\frac{2u + 4v}{u+v}\right) = \left(\frac{2v}{u+v}\right) = \left(\frac{2}{u+v}\right)^{e_v+1}\left(\frac{v'}{u+v}\right).
$$
Using that $u\equiv 1\bmod 8$ and $v\equiv 0\bmod 4$, so that $\left(\frac{2}{u+v}\right) = \left(\frac{2}{1+v}\right)$ and $\left(\frac{v'}{u+v}\right)=\left(\frac{v'}{u}\right)=\left(\frac{v}{u}\right)$, we arrive at 
$$
\left(\frac{2h}{g}\right) = \left(\frac{2}{1+v}\right)^{e_v+1}\left(\frac{v}{u}\right).
$$
As
$$
\left(\frac{2}{1+v}\right)^{e_v+1} = 
\begin{cases}
1 & \text{if }v\equiv 0 \bmod 8 \\
-1 & \text{if }v\equiv 4 \bmod 8,
\end{cases}
$$
we obtain the formula
\begin{equation}\label{lastfactor}
\left(\frac{2h}{g}\right) = (-1)^{v/4}\left(\frac{v}{u}\right).
\end{equation}

\subsection{The formula for $\beta_p$}
Recall from \eqref{defapbp} that $\beta_p$ is defined to be
$$
\beta_p = (-1)^{(a-1+b+2d+h(-p)+h(-2p))/8}.
$$
Combining \eqref{overZ8}, \eqref{firstfactor}, \eqref{middlefactor}, and \eqref{lastfactor}, we arrive at the formula
$$
\beta_p = (-1)^{(a-1)/8}(-1)^{v/4}\left(\frac{v}{u}\right).
$$
Using \eqref{z1} and \eqref{z2}, namely that $p = a^2+b^2 \equiv 1\bmod 16$ with $a-1\equiv b\equiv 0\bmod 8$, one can check that
$$
\frac{a-1}{8}\equiv \frac{p-1}{16} \bmod 2.
$$
We thus arrive at the final form of the formula for $\beta_p$ that we will subsequently extend into an oscillating sequence to be sieved for primes:
\begin{equation}\label{finalbp}
\beta_p = (-1)^{(p-1)/16}(-1)^{v/4}\left(\frac{v}{u}\right).
\end{equation}

\subsection{A comment on the formula for $\beta_p$}\label{catastrophe}
Due to the presence of fourth power residue symbols in the algebraic criteria of Leonard and Williams, it seems most natural to define $\beta_p$ over $\QQ(\zeta_8)$ as in \eqref{overZ8}. However, the final formula \eqref{finalbp} suggests that $\beta_p$ also has a natural definition over the smaller field $\QQT$, which is very advantageous. Roughly speaking, if we were only able to define $\beta_p$ over $\QQ(\zeta_8)$, our analytic arguments would require us to prove significant cancellation (namely power-saving in $X$) in sums resembling 
$$
\sum_{\substack{\alpha = a_0 + a_1\zeta_8+a_2\zeta_8^2 + a_3\zeta_8^3\\ a_i\in\ZZ,\ |a_i|\leq X^{1/4} \\ \alpha\sigma(\alpha) = u+v\sqrt{2},\ u\text{ odd}, >0}}\left(\frac{v}{u}\right).
$$
This appears to be a very difficult problem that we see how to solve only by appealing to a standard conjecture on short character sums that is just out of reach of the deep Burgess's inequality (see Conjecture $C_n$ for $n=4$ from \cite{FIMR}). Instead, working over $\QQ(\sqrt{2})$ leads to sums of the form
$$
\sum_{\substack{u, v\in \ZZ,\ u\text{ odd},>0\\ 2|v| < u \leq X^{1/2}}}\left(\frac{v}{u}\right),
$$
which can be handled even with just the classical P\'{o}lya-Vinogradov inequality.

Moreover, we remark that although the essential spin factor in $\beta_p$ comes from the formula for $(-1)^{(h(-p)+h(-2p))/8}$ (see the last factor in \eqref{overZ8}), the twist by $(-1)^{(a-1+b+2d)/8}$ in the definition of $\beta_p$ (see \eqref{defapbp}) is exactly what allows us to cancel the seemingly innocuous but fatal first factor appearing in \eqref{overZ8}. In fact, determining the symbol
$$
\left(\frac{2\sqrt{2}\ve^{-1}}{\varpi}\right)_{K, 4} = \left(\frac{2(2-\sqrt{2})}{\varpi}\right)_{K, 4}
$$
is tantamount to determining the splitting of $p$ in the extension $L/\QQ$, with
$$
L = \QQ\left(\zeta_8, \sqrt[4]{2-\sqrt{2}}\right).
$$
The key observation is that the Galois group of $L/\QQT$ is isomorphic to the quaternion group $Q_8$, which is \textit{not abelian}; hence, by class field theory, the first factor in \eqref{overZ8} \textit{cannot} be determined by congruence conditions on a prime lying above $p$ in $\ZZT$. This would then force us to work over the bigger field $\QQ(\zeta_8)$. We challenge the reader to give an unconditional proof that there exists a $\delta>0$ such that for all $X\geq 3$ we have
$$
\sum^{\ast}_{p\leq X}(-1)^{(h(-p)+h(-2p))/8}\ll X^{1-\delta},
$$  
where, as before, $\ast$ restricts the summation to primes $p$ that split completely in $\QQ(\zeta_{16}, \sqrt[4]{2})$.

\section{Construction of the spin sequence}
Following the method in \cite{Milovic2}, we wish to construct a bounded sequence of complex numbers $\{b_{\nn}\}_{\nn}$ indexed by non-zero ideals $\nn$ of $\ZZT$ such that $b_{\pp} = \beta_p$ whenever $\pp$ is a prime ideal lying above a prime $p$ that splits completely in $\QQ(\zeta_{16}, \sqrt[4]{2})/\QQ$ and such that we can prove power-saving estimates for sums of the form
$$
\sum_{\substack{\Norm\nn\leq X\\ \nn\equiv 0\bmod \dd}}b_{\nn} \quad\quad \text{and}\quad\quad\sum_{\Norm\mm\leq M}\sum_{\Norm\nn\leq N}v_{\mm}w_{\nn}b_{\mm\nn}.
$$
To this end, we will now state the key result from \cite{Milovic2} that lets us do so. Let $\ve = 1+\sqrt{2}$, as before. For an \textit{odd, totally positive} (and not necessarily prime) element $\alpha = u+v\sqrt{2}\in\ZZT$, define
$$
[\alpha] = \left(\frac{v}{u}\right).
$$
Then \cite[Proposition 2, p.\ 979]{Milovic2} implies that
\begin{equation}\label{key1}
[\alpha] = [\ve^8\alpha].
\end{equation}
Now, since $\ZZT$ is a principal ideal domain and $\ve$ is a unit of norm $-1$, every non-zero ideal of $\ZZT$ can be generated by a totally positive element. Let $\alpha = u+v\sqrt{2}$ be a totally positive generator of a non-zero ideal $\nn$; all the other totally positive generators of $\nn$ are of the form $\ve^{2k}\alpha$ for some integer $k$. Suppose now that $\nn$ is odd, i.e., that $\Norm\nn$ is odd. Then \eqref{key1} implies that the quantity 
$$
[\alpha]+[\ve^2\alpha]+[\ve^4\alpha]+[\ve^6\alpha]
$$
depends only on the ideal $\nn = \alpha\ZZT$ and \textit{not} on the choice of the totally positive generator $\alpha$ of $\nn$. Eventually we will need to be able to detect $\alpha$ with $u$ and $v$ satisfying $u-1\equiv 2v\equiv 0 \bmod 8$ (see \eqref{z2}) as well as detect when $v\equiv 0\bmod 8$ (to control the factor $(-1)^{v/4}$ in \eqref{finalbp}). The idea is to first detect $\alpha$ with $\Norm(\alpha) =u^2-2v^2\equiv 1\bmod 16$ via multiplicative Dirichlet characters modulo $16$. Then we will detect when $u\equiv 1 \bmod 8$; this already ensures that $v\equiv 0\bmod 4$ provided that $u^2-2v^2\equiv 1\bmod 16$. To detect when $u\equiv 1\bmod 8$, we should study how $(u\bmod 8, v\bmod 4)$ changes as we multiply $\alpha$ by successive powers of $\ve^2$. We compute that
$$
\ve^2\alpha = (3+2\sqrt{2})(u+v\sqrt{2}) = (3u+4v)+(2u+3v)\sqrt{2}.
$$
The orbits of the map $(u \bmod 8, v\bmod 4)\mapsto (3u+4v \bmod 8, 2u+3v \bmod 4)$ for $u$ odd and $v$ even can be listed as follows (note that $u$ odd implies that $2u\equiv 2\bmod 4$):
\begin{align*}
& (u, 0)\mapsto (3u, 2)\mapsto (u, 0) \\
& (u, 2)\mapsto (3u, 0)\mapsto (u, 2).
\end{align*}
Hence if $u'+v'\sqrt{2} = \ve^2(u+v\sqrt{2})$ with $v, v'$ even, then either $\{u\bmod 8, u'\bmod 8\} = \{1\bmod 8, 3\bmod 8\}$ or $\{u\bmod 8, u'\bmod 8\} = \{5\bmod 8, 7\bmod 8\}$.

Now, for the rest of the paper, fix a square root of $-1$ and denote it by $i$. We define, for each Dirichlet character $\chi$ modulo $8$ and each odd totally positive element $u+v\sqrt{2}\in\ZZT$ with $v$ even,
$$
[u+v\sqrt{2}]_{\chi} = i^{v/2}\chi(u)\left(\frac{v}{u}\right).
$$
Next, for each pair of Dirichlet characters $\chi$ modulo $8$ and $\psi$ modulo $16$ and each ideal $\nn$ satisfying $\Norm\nn \equiv 1\bmod 8$, we define
\begin{equation}\label{spin}
b_{\nn}(\chi, \psi) = \frac{1}{2}\psi(\Norm(\nn))i^{(\Norm(\nn)-1)/8}\left([\alpha]_{\chi}+[\ve^2\alpha]_{\chi}+[\ve^4\alpha]_{\chi}+[\ve^6\alpha]_{\chi}\right),
\end{equation}
where $\alpha$ is any totally positive generator of $\nn$. To prove that the right-hand side above is a well-defined function of $\nn$, it suffices to show that $[\ve^8\alpha]_{\chi} = [\alpha]_{\chi}$. Indeed, $\ve^8\alpha = (577+408\sqrt{2})(u+v\sqrt{2}) = (577u+816v)+(408u+577v)\sqrt{2}$, so
$$
i^{(408u+577v)/2} = i^{v/2}\cdot i^{4(51u+72v)} = i^{v/2},\quad \chi(577u+816v) = \chi(u+8(72u+102v)) = \chi(u),
$$
and hence
$$
[\ve^8\alpha]_{\chi} =  i^{(408u+577v)/2}\chi(577u+408v)[\ve^8\alpha] = i^{v/2}\chi(u)[\alpha] = [\alpha]_{\chi}.
$$
We extend $b_{\nn}(\chi, \psi)$ to the remaining non-zero ideals by $0$, i.e., we set  $b_{\nn}(\chi, \psi)=0$ whenever $\Norm\nn\not\equiv 1\bmod 8$.

Finally, when $\nn = \pp$ is a prime ideal in $\ZZT$ lying above a prime $p\equiv \pm 1\bmod 8$, i.e., above a prime that splits completely in $\ZZT$, then
\begin{equation}\label{agreement}
\frac{1}{4\cdot 8}\sum_{\chi\bmod 8}\sum_{\psi\bmod 16}b_{\pp}(\chi, \psi) = 
\begin{cases}
\beta_p & \text{if }p\text{ splits completely in }\QQ(\zeta_{16}, \sqrt[4]{2}) \\
0 & \text{otherwise.}
\end{cases}
\end{equation}

\section{Proof of Proposition~\ref{PropBp}}
To prove Proposition~\ref{PropBp}, it suffices to prove that 
\begin{equation}\label{reductionBp}
\sum_{\Norm\pp\leq X}b_{\pp}(\chi, \psi) \ll X^{1-\frac{1}{200}}
\end{equation}
for all pairs of Dirichlet characters $\chi$ modulo $8$ and $\psi$ modulo $16$. Indeed, if this is the case, since the contribution from $\sqrt{2}\ZZT$ and the inert primes is $\ll X^{\frac{1}{2}}$ when we order the prime ideals by norm, we see that the same estimate holds when the sum is restricted to split primes $\pp$. Then adding together finitely many such sums, one for each pair of Dirichlet characters $\chi$ modulo $8$ and $\psi$ modulo $16$, we obtain the sum in Proposition~\ref{PropBp} via formula \eqref{agreement}.
The proof of \eqref{reductionBp} is essentially no different from the proof of \cite[Theorem 3, p.994]{Milovic2}. In short, one uses Vinogradov's method of sums of type I and type II, in the form presented in \cite[Section 5]{FIMR}. One first selects a suitable fundamental domain in $\ZZ^2$ for the action of $\ve^2$ on non-zero elements of $\ZZT$, so that each point in the domain corresponds to exactly one non-zero ideal in $\ZZT$. This allows us to pass to sums over \textit{elements} $\alpha$ in the domain (or a translate thereof) instead of \textit{ideals}; the summand then becomes $(-1)^{\Norm((\alpha)-1)/8}\psi(\Norm(\alpha))[\alpha]_{\chi}$. To control the factors appearing in the summand other than $[\alpha]$, one breaks up the relevant sum according to the congruence class of $\alpha$ modulo $16$, shows the desired estimate for each such sum, and then adds together the contributions from the finitely many congruence classes. In this way, one proves that 
$$
\sum_{\substack{\Norm\nn\leq X\\ \nn\equiv 0\bmod \dd}}b_{\nn}(\chi, \psi) \ll_{\epsilon} X^{\frac{5}{6}+\epsilon}
$$
uniformly in $\dd$. The key idea is to breaking up the domain into horizontal segments and apply the P\'{o}lya-Vinogradov inequality to the sum over each segment.  See \cite[Section 5]{Milovic2} for details. Next, one proves that
$$
\sum_{\Norm\mm\leq M}\sum_{\Norm\nn\leq N}v_{\mm}w_{\nn}b_{\mm\nn}(\chi, \psi) \ll_{\epsilon}(MN)^{\frac{11}{12}+\epsilon}\left(M^{\frac{1}{12}}+N^{\frac{1}{12}}\right)
$$ 
uniformly for all bounded sequences of complex numbers $\{v_{\mm}\}$ and $\{w_{\nn}\}$ indexed by non-zero ideals of $\ZZT$. One again isolates the key spin $[\alpha\beta]$ from the other factors in the definition of $[\alpha\beta]_{\chi}$ by restricting the congruence classes of $\alpha$ and $\beta$ modulo $16$, and then uses the very key \cite[Proposition 8, p.1010]{Milovic2} (modeled after \cite[Lemma 20.1, p.1021]{FI1}) to factor $[\alpha\beta]$ into essentially $[\alpha]$, $[\beta]$, which are absorbed into the sequences $v_{\mm}$ and $w_{\nn}$, and a quadratic residue symbol $\left(\alpha/\beta\right)_{\QQT, 2}$. One finishes by applying the double oscillation result \cite[Lemma 22, p.1009]{Milovic2}. See \cite[Section 6]{Milovic2} for details. The result then follows by applying \cite[Proposition 5.2, p.722]{FIMR}.

\section{Proof of Propositions~\ref{PropAp} and \ref{PropApBp}}
Finally, we explain how to obtain Propositions~\ref{PropAp} and \ref{PropApBp} from \cite{KM1} and \cite{Koymans}, respectively.
\subsection{Proof of Proposition~\ref{PropAp}}
In \cite{KM1}, it is proved that 
$$
\sum^{'}_{p\leq X}\alpha_p\ll X^{1-\frac{1}{3200}},
$$
where $'$ restricts the summation to primes $p$ that split completely in $\QQ(\zeta_8, \sqrt[4]{2})$. The proof has a similar structure to the proof in \cite{Milovic2}, but one works over $\ZZ[\zeta_8]$ instead of $\ZZT$ (although we said in Section~\ref{catastrophe} that working over $\ZZ[\zeta_8]$ would only give a conditional result as it pertains to Proposition~\ref{PropBp}, there are other circumstances we were able to exploit in \cite{KM1} to obtain an unconditional result). So if we choose an element $\varpi\in\ZZ[\zeta_8]$ whose norm is a prime $p$ that splits completely in $\QQ(\zeta_8, \sqrt[4]{2})$, then $p$ further splits in $\QQ(\zeta_{16}, \sqrt[4]{2})$ if and only if the quadratic residue symbol $(\zeta_8/\varpi)_{\QQ(\zeta_8), 2}$ is equal to $1$ (since $\zeta_{16}$ is a square root of $\zeta_8$). Hence we simply multiply the symbol $a(\chi)_{\nn}$ in \cite[(2.3), p.6]{KM1} by the indicator function 
$$
\frac{1}{2}\left(1+\left(\frac{\zeta_8}{\nn}\right)_{\QQ(\zeta_8), 2}\right).
$$
After expanding, we see that one also has to prove \cite[Proposition 3.7, p. 13]{KM1} and \cite[Proposition 3.8, p. 14]{KM1} with $a(\chi)_{\nn}$ replaced by $a(\chi)_{\nn}(\zeta_8/\nn)_{\QQ(\zeta_8), 2}$. By quadratic reciprocity, the quadratic residue symbol $(\zeta_8/\nn)_{\QQ(\zeta_8), 2}$ is controlled by the congruence class modulo $8\ZZ[\zeta_8]$ of a generator $\alpha$ of $\nn$, and so the same proofs apply (since in the proofs in \cite{KM1} one immediately reduces to sums over a fixed congruence class modulo $F = 16$). 

\subsection{Proof of Proposition~\ref{PropApBp}}
In \cite{KM1}, it is proved that 
$$
\sum^{\ast}_{p\leq X}\alpha_p\beta_p\ll X^{1-\delta},
$$
where $\ast$ restricts the summation to primes $p$ that split completely in $\QQ(\zeta_{16}, \sqrt[4]{2})$, but the result is conditional on Conjecture $C_n$ for $n=8$ from \cite{FIMR}. The reason we needed ths conjecture is that this time we carried out the analytic estimates over $\QQ(\zeta_8, \sqrt{1+\zeta_4})$, a number field of degree $8$ over $\QQ$. In proving a result on $(-1)^{(h(-p))/8}$, Koymans \cite{Koymans} managed to work over $\QQ(\zeta_8)$ (where he again encountered extenuating circumstances that allowed him to obtain an unconditional result, without having to use Conjecture $C_n$ for $n = 4$). In \cite[Section 7]{KM2}, we explained how to account for the twist by $(-1)^{(a-1+b+2d)/8}$ (see \eqref{defapbp}). When working over the smaller field $\QQ(\zeta_8)$, this is even easier. As shown above in Section~\ref{SectionAC}, the twist is actually equal to 
$$
(-1)^{\frac{p-1}{16}}\left(\frac{2-\sqrt{2}}{\varpi}\right)_{\QQ(\zeta_8), 4}
$$
and hence controlled by the congruence class of $\varpi$ modulo $32\ZZ[\zeta_8]$ ($(-1)^{(p-1)/16}$ is determined by $p\bmod 32$, which is certainly determined by $\varpi\bmod 32\ZZ[\zeta_8]$; the factor $(2-\sqrt{2}/\varpi)_{\QQ(\zeta_8), 4}$ is determined by $\varpi\bmod 8\ZZ[\zeta_8]$, by quartic reciprocity). As Koymans already restricts the sums appearing in \cite{Koymans} to fixed congruence classes modulo a much higher power of $2$, his proof safely carries through also for $(-1)^{(h(-p))/8}$ twisted by $(-1)^{(a-1+b+2d)/8}$.

\bibliographystyle{abbrv}
\bibliography{M4_References}

\begin{thebibliography}{10}

\bibitem{FI1}
J.~B. Friedlander and H.~Iwaniec.
\newblock The polynomial {$X^2+Y^4$} captures its primes.
\newblock {\em Ann. of Math. (2)}, 148(3):945--1040, 1998.

\bibitem{FIMR}
J.~B. Friedlander, H.~Iwaniec, B.~Mazur, and K.~Rubin.
\newblock The spin of prime ideals.
\newblock {\em Invent. Math.}, 193(3):697--749, 2013.

\bibitem{KW84}
P.~Kaplan and K.~S. Williams.
\newblock On the strict class number of {${\bf Q}(\sqrt{2p})$} modulo {$16,$}
  {$p\equiv 1$} {$({\rm mod}\,8)$} prime.
\newblock {\em Osaka J. Math.}, 21(1):23--29, 1984.

\bibitem{Koymans}
P.~{Koymans}.
\newblock {The $16$-rank of $\mathbb{Q}(\sqrt{-p})$}.
\newblock {\em ArXiv e-prints}, page arXiv:1809.07167, Sept. 2018.

\bibitem{KM1}
P.~Koymans and D.~Z. Milovic.
\newblock On the 16-rank of class groups of $\mathbb{Q}(\sqrt{-2p})$ for primes
  $p \equiv 1 \bmod 4$.
\newblock {\em International Mathematics Research Notices}, page rny010, 2018.

\bibitem{KM2}
P.~Koymans and D.~Z. Milovic.
\newblock Spins of prime ideals and the negative pell equation
  $x^{2}-2py^{2}=-1$.
\newblock {\em Compositio Mathematica}, 155(1):100--125, 2019.

\bibitem{Lemm}
F.~Lemmermeyer.
\newblock Ideal class groups of cyclotomic number fields. {I}.
\newblock {\em Acta Arith.}, 72(4):347--359, 1995.

\bibitem{LW82}
P.~A. Leonard and K.~S. Williams.
\newblock On the divisibility of the class numbers of {$Q(\sqrt{-p})$} and
  {$Q(\sqrt{-2p})$} by {$16$}.
\newblock {\em Canad. Math. Bull.}, 25(2):200--206, 1982.

\bibitem{Milovic2}
D.~Z. Milovic.
\newblock On the 16-rank of class groups of {$\Bbb Q(\sqrt{-8p})$} for
  {$p\equiv-1 \mod4$}.
\newblock {\em Geom. Funct. Anal.}, 27(4):973--1016, 2017.

\bibitem{Redei2}
L.~R{\'e}dei.
\newblock Ein neues zahlentheoretisches {S}ymbol mit {A}nwendungen auf die
  {T}heorie der quadratischen {Z}ahlk\"orper. {I}.
\newblock {\em J. Reine Angew. Math.}, 180:1--43, 1939.

\bibitem{Scholz}
A.~Scholz.
\newblock \"{U}ber die {L}\"osbarkeit der {G}leichung {$t^2-Du^2=-4$}.
\newblock {\em Math. Z.}, 39(1):95--111, 1935.

\bibitem{Ste1}
P.~Stevenhagen.
\newblock Ray class groups and governing fields.
\newblock In {\em Th\'eorie des nombres, {A}nn\'ee 1988/89, {F}asc.\ 1}, Publ.
  Math. Fac. Sci. Besan\c con, page~93. Univ. Franche-Comt\'e, Besan\c con,
  1989.

\bibitem{Ste2}
P.~Stevenhagen.
\newblock Divisibility by {$2$}-powers of certain quadratic class numbers.
\newblock {\em J. Number Theory}, 43(1):1--19, 1993.

\bibitem{SteConj}
P.~Stevenhagen.
\newblock The number of real quadratic fields having units of negative norm.
\newblock {\em Experiment. Math.}, 2(2):121--136, 1993.

\end{thebibliography}

\end{document}